\documentclass[nosumlimits,twoside]{amsart}%
\usepackage{amssymb}
\usepackage{amsmath}
\usepackage[bookmarksnumbered,plainpages,hyperindex]{hyperref}
\usepackage{amsfonts}
\usepackage{float}
\usepackage{graphicx}
\newtheorem{theorem}{\sc Theorem}[section]
\newtheorem{proposition}[theorem]{\sc Proposition}
\newtheorem{lemma}[theorem]{\sc Lemma}
\newtheorem{corollary}[theorem]{\sc Corollary}
\theoremstyle{definition}
\newtheorem{definition}[theorem]{\sc Definition}
\newtheorem{definitions}[theorem]{\sc Definitions}

\theoremstyle{remark}

\newtheorem{claim}[theorem]{}
\newtheorem{noname}[theorem]{}

\def\tl{{\vartriangleleft}}
\def\tr{{\vartriangleright}}

\def\ot{\otimes}

\input{xy}
\xyoption{all}
\newcommand{\diagBradedl}{\xymatrix@R=30pt@C=3pt{
&B\otimes B\ot B\otimes B \ar[rr]^{ B\ot c\ot B} &&B\otimes B\ot
B\otimes B
\ar[dr]|{ m\ot m}&\\
B\ot B\ar[ur]|{\Delta\ot \Delta } \ar[drr]|{m} &&&&B\otimes B
\\&&B\ar[urr]|{\Delta }&&&}}
\newcommand{\diagBradedr}{\xymatrix@C=40pt{
B\otimes B \ar[d]_{m} \ar[r]^{\varepsilon \otimes \varepsilon} & \mathbf{1}\ot \mathbf{1} \ar[d]^{m_{\mathbf{1}}} \\
B \ar[r]^{\varepsilon} & \mathbf{1}   }}
\newcommand{\diagCosepl}{\xymatrix@R=15pt@C=50pt{
{^C\mathcal{M}^C} \ar[d]_{\mathbb{T}} \ar[r]^{F'} & {^{C'}\mathcal{M'}^{C'}} \ar[d]^{\mathbb{T}'} \\
\mathcal{M}  \ar[r]_{F} & \mathcal{M'}   }}
\newcommand{\diagCosepr}{\xymatrix@R=15pt@C=50pt{
{^C\mathcal{M}^C} \ar[r]^{F'} & {^{C'}\mathcal{M'}^{C'}}  \\
\mathcal{M}  \ar[u]^{\mathbb{H}} \ar[r]_{F} & \mathcal{M'}
\ar[u]_{\mathbb{H}'}  }}
\newcommand{\diagFSl}{\xymatrix@R=15pt@C=50pt{
{^B\mathfrak{M}_B^B} \ar[d]_{\mathbb{T}} \ar[r]^{F'} & {^{B}\mathfrak{M}^{B}} \ar[d]^{\mathbb{T}'} \\
\mathfrak{M}_B  \ar[r]_{F} & \mathfrak{M}   }}
\newcommand{\diagFSr}{\xymatrix@R=15pt@C=50pt{
{^B\mathfrak{M}^B_B} \ar[r]^{F'} & {^{B}\mathfrak{M}^{B}}  \\
\mathfrak{M}_B  \ar[u]^{\mathbb{H}} \ar[r]_{F} & \mathfrak{M}
\ar[u]_{\mathbb{H}'}  }}
\begin{document}
\title{Weak Projections onto a Braided Hopf Algebra}
\author{A. Ardizzoni}
\address{University of Ferrara, Department of Mathematics, Via Machiavelli 35, Ferrara,
I-44100, Italy}
\email{rdzlsn@unife.it}
\urladdr{http://www.unife.it/utenti/alessandro.ardizzoni}
\author{C. Menini}
\address{University of Ferrara, Department of Mathematics, Via Machiavelli 35, Ferrara,
I-44100, Italy}
\email{men@unife.it}
\urladdr{http://www.unife.it/utenti/claudia.menini}
\author{D. \c{S}tefan}
\address{University of Bucharest, Faculty of Mathematics, Str. Academiei 14, Bucharest,
Ro--70109, Romania}
\email{dstefan@al.math.unibuc.ro}
\thanks{This paper was written while the first two authors were members of G.N.S.A.G.A. with
partial financial support from M.I.U.R..}
\thanks{The third author was financially supported by the project
CeEx 2006.}

\begin{abstract}
We show that, under some mild conditions, a bialgebra in an
abelian and coabelian braided monoidal category has a weak
projection onto a formally smooth (as a coalgebra) sub-bialgebra
with antipode; see Theorem \ref{teo: Magnum}. In the second part
of the paper we prove that bialgebras with weak projections are
cross product bialgebras; see Theorem \ref{te: bialgebras}. In the
particular case when the bialgebra $A$ is cocommutative and a
certain cocycle associated to the weak projection is trivial we
prove that $A$ is a double cross product, or biproduct in Madjid's
terminology. The last result is based on a universal property of
double cross products which, by Theorem \ref{te: matched pair},
works in braided monoidal categories. We also investigate the
situation when the right action of the associated matched pair is
trivial.

\end{abstract}
 \maketitle

\section*{Introduction}
Hopf algebras in a braided monoidal category are very important
structures. Probably the first known examples are
$\mathbb{Z}$-graded and $\mathbb{Z}_2$-graded bialgebras (also
called superbialgebras), that already appeared in the work of
Milnor-Moore and MacLane. Other examples, such as bialgebras in
the category of Yetter-Drinfeld modules, arose in a natural way in
the characterization as a double crossed product of (ordinary)
Hopf algebras with a projection \cite{Rad} . Some braided
bialgebras have also played a central role in the theory of
quantum groups.

The abundance of examples and their applications explain the
increasing interest for these objects and the attempts in
describing their structure. For example in \cite{BD, BD-Cross1,
BD-Cross2, Schauenburg1} several generalized versions of the
double cross product  bialgebra in a braided monoidal category
$\mathfrak{M}$, generically called cross product bialgebras, are
constructed. All of them have the common feature that, as objects
in $\mathfrak{M}$, they are the tensor product of two objects in
$\mathfrak{M}$. Let $A$ be such a cross product, and let $R$ and
$B$ the corresponding objects such that $A\simeq R\ot B$.
Depending on the particular type of cross product, the objects $R$
and  $B$ may have additional properties, like being algebras
and/or coalgebras. These structures may also satisfy some
compatibility relations. For example, we can look for those cross
product bialgebras $A\simeq R\ot B$ such that there are a
bialgebra morphism $\sigma:B\to A$ and a right $B$-linear
coalgebra map $\pi:A\to B$ satisfying the relation
$\pi\sigma=\mathrm{Id}_B$ (here $A$ is a $B$-module via $\sigma$).
For simplicity, we will say that $A$ is a bialgebra (in
$\mathfrak{M}$) with weak projection $\pi$ on $B$. In the case
when $\mathfrak{M}$ is the category of vector spaces, the problem
of characterizing bialgebras with a weak projection was considered
in \cite{Schauenburg1}.

The purpose of this is two fold. We assume that $\mathfrak{M}$ is
a semisimple abelian and coabelian braided monoidal category (see
Definition~\ref{abelianmonoidal}) and that $\sigma:B\to A$ is
morphism of bialgebras in $\mathfrak{M}$. First, we want to show
that there is a retraction $\pi$ of $\sigma$ which is a right
$B$-linear morphism of coalgebras, provided that $B$ is formally
smooth as a coalgebra and that the $B$-adic coalgebra filtration
on $A$ is exhaustive (see Theorem~\ref{teo: Magnum}). Secondly,
assuming that a retraction $\pi$ as above exists, we want to show
that $A$ is factorizable and to describe the corresponding
structure $R$ that arises in this situation. For this part of the
paper we use the results of \cite{BD-Cross2}, that help us to
prove that $A$ is the cross product algebra $R\rtimes B$ where $R$
is the `coinvariant' subobject with respect to the right coaction
of $B$ on $A$ defined by $\pi$ (see Theorem~\ref{te: bialgebras}).
Several particular cases are also investigated. In
Theorem~\ref{te: double cross product}, under the additional
assumption that $A$ is cocommutative and a certain cocyle is
trivial, we describe the structure of $A$ as a biproduct bialgebra
of a certain matched pair (see Theorem~\ref{te: matched pair}).

Finally, we would like to note that some applications of the last
mentioned result and its corollaries (see Proposition~\ref{pr:
bosonization}) are given in \cite{AMS3}. As a matter of fact, our
interest for the problems that we study in this paper originates
in our work on the structure of cocommutative Hopf algebras with
dual Chevalley property from \cite{AMS3}. In particular,
\cite[Theorem 6.14]{AMS3} and \cite[Theorem 6.16]{AMS3} are direct
consequences of the main results of this article.

\section{Hopf algebras in a braided category $\mathfrak{M}$}

\begin{definition}
\label{abelianmonoidal}An \emph{abelian monoidal category} is a monoidal
category $(\mathfrak{M},\otimes,\mathbf{1})$ such that:

\begin{enumerate}
\item $\mathcal{M}$ is an abelian category

\item both the functors $X\otimes(-):\mathfrak{M}\rightarrow\mathfrak{M}$ and
$(-)\otimes X:\mathfrak{M}\rightarrow\mathfrak{M}$ are additive and right
exact, for every object $X\in\mathfrak{M}$.
\end{enumerate}

A \emph{coabelian monoidal category} is a monoidal category $(\mathfrak{M}%
,\otimes,\mathbf{1})$ such that:

\begin{enumerate}
\item $\mathcal{M}$ is an abelian category

\item both the functors $X\otimes(-):\mathfrak{M}\rightarrow\mathfrak{M}$ and
$(-)\otimes X:\mathfrak{M}\rightarrow\mathfrak{M}$ are additive and left
exact, for every object $X\in\mathfrak{M}$.
\end{enumerate}
\end{definition}

\begin{definitions}
\label{def braiding} A \emph{braided monoidal category} $(\mathfrak{M}%
,\otimes,\mathbf{1},c)$ is a monoidal category $(\mathfrak{M}%
,\otimes,\mathbf{1})$ equipped with a \emph{braiding} $c$, that is a natural
isomorphism
\[
c_{X,Y}:X\otimes Y\rightarrow Y\otimes X
\]
satisfying
\[
c_{X\otimes Y,Z}=(c_{X,Z}\otimes Y)(X\otimes c_{Y,Z})\qquad\text{and}\qquad
c_{X,Y\otimes Z}=(Y\otimes c_{X,Z})(c_{X,Y}\otimes Z).
\]
For further details on these topics, we refer to \cite[Chapter XIII]{Ka}.

\label{cl: brdBialg} A \emph{bialgebra} $(B,m,u,\Delta,\varepsilon)$ in a
braided monoidal category $(\mathfrak{M},{\otimes} ,\mathbf{1},c)$ consists of
an algebra $(B,m,u)$ and a coalgebra $(B,\Delta,\varepsilon)$ in
$\mathfrak{M}$ such that the diagrams in Figure~\ref{fig:bialgebra} are
commutative. \begin{figure}[h]
\begin{center}
\fbox{\includegraphics{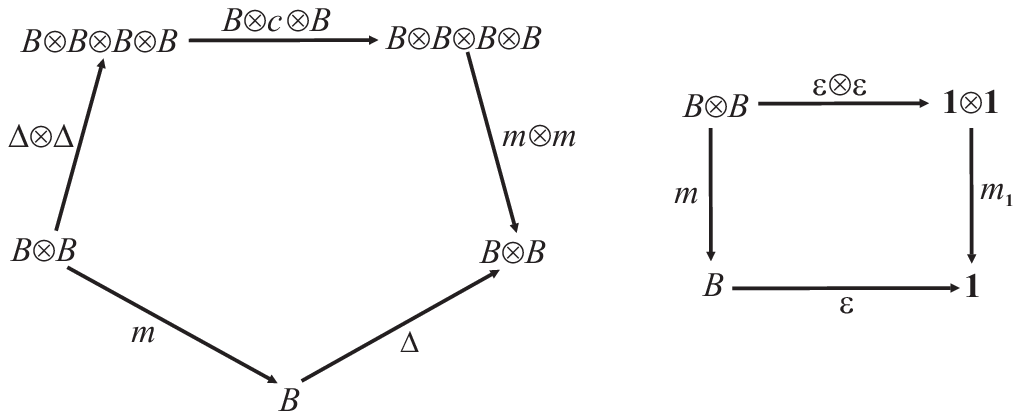}}
\end{center}
\caption{The definition of bialgebras in $\mathfrak{M}$.}\label{fig:bialgebra}%
\end{figure}
\end{definitions}

\begin{noname}
For any bialgebra $B$ in a monoidal category $(\mathfrak{M},\otimes
,\mathbf{1})$ we define the monoidal category $\left(  \mathfrak{M}%
_{B},\otimes,\mathbf{1}\right)  $ of right $B$-modules in $\mathfrak{M}$ as in
\cite{BD}. The tensor product of two right $B$-modules $\left(  M,\mu
_{M}\right)  $ and $\left(  N,\mu_{N}\right)  $ carries a right $B$-module
structure defined by:%
\[
\mu_{M\otimes N}=(\mu_{M}\otimes\mu_{N})(M\otimes c_{N,B}\otimes B)(M\otimes
N\otimes\Delta).
\]
Moreover, if $(\mathfrak{M},\otimes,\mathbf{1})$ is an abelian
monoidal category, then $\left(
\mathfrak{M}_{B},\otimes,\mathbf{1}\right)  $ is an abelian
monoidal category too. Assuming that $(\mathfrak{M},\otimes
,\mathbf{1})$ is abelian and coabelian one proves that $\left(
\mathfrak{M}_{B},\otimes ,\mathbf{1}\right) $ is coabelian too.

Obviously, $(B,\Delta,\varepsilon) $ is a coalgebra both in $\left(
\mathfrak{M}_{B},\otimes,\mathbf{1}\right)  $ and $\left(  _{B}\mathfrak{M}%
_{B},\otimes,\mathbf{1}\right)  $. Of course, in both cases, $B$ is regarded
as a left and a right $B$-module via the multiplication on $B$.
\end{noname}

\begin{claim}
To each coalgebra $(C,\Delta,\varepsilon)$ in $(\mathcal{M},\otimes
,\mathbf{1},a,l,r)$ one associates a class of monomorphisms
\[
{^{C}\mathcal{I}^{C}} :=\{g\in{^{C}\mathcal{M}^{C}}\mid\exists f\text{ in
}\mathcal{M}\text{ s.t. } fg=\mathrm{Id}\}.
\]
Recall that $C$ is \emph{coseparable} whenever the comultiplication $\Delta$
cosplits in $^{C}\mathcal{M}^{C}$. We say that $C$ is \emph{formally smooth
}in ${\mathcal{M}}$ if $\text{Coker}\Delta_{C}$ is ${^{C}\mathcal{I}^{C}}%
$-injective. For other characterizations and properties of coseparable and
formally smooth coalgebras the reader is referred to \cite{AMS1} and
\cite{Ar1}. In the same papers one can find different equivalent definitions
of separable functors.
\end{claim}

\begin{claim}
\label{ex: co monoidal functor} Let $(F,\phi_{0},\phi_{2}):(\mathcal{M}%
,\otimes,\mathbf{1},a,l,r\mathbf{)\rightarrow(}\mathcal{M}^{\prime}%
,\otimes^{\prime},\mathbf{1}^{\prime},a^{\prime},l^{\prime},r^{\prime
}\mathbf{)}$ be a monoidal functor between two monoidal categories, where
$\phi_{2}(U,V):F(U\otimes V)\rightarrow F(U)\otimes^{\prime}F(V),$ for any
$U,V\in\mathcal{M}$ and $\phi_{0}:\mathbf{1}^{\prime}\mathbf{\rightarrow
F(1)}.$ If $(C,\Delta,\varepsilon)$ is a coalgebra in $\mathcal{M}$ then
$(C^{\prime},\Delta_{C^{\prime}},\varepsilon_{C^{\prime}}):=(F(C),\Delta
_{F(C)},\varepsilon_{F(C)})$ is a coalgebra in $\mathcal{M}^{\prime}$, with
respect to the comultiplication and the counit given by
\[
\Delta_{F(C)}:=\phi_{2}^{-1}(C,C)F(\Delta), \qquad\varepsilon_{F(C)}:=\phi
_{0}^{-1}F(\varepsilon).
\]
Let us consider the functor $F^{\prime}:{}^{C}\!\mathcal{M}^{C}\rightarrow
{^{C^{\prime}}\!\mathcal{M}^{\prime}{}^{C^{\prime}}}$ that associates to
$(M,{^{C}\!\rho_{M}},{\rho_{M}^{C}})$ the object $(F(M),{^{C^{\prime}}%
\!\!\rho_{F(M)}},{\rho_{F(M)}^{C^{\prime}}}),$ where
\[
{^{C^{\prime}}\!\!\rho_{F(M)}}:=\phi_{2}^{-1}(C,M)F(\,^{C}\!\rho_{M}%
),\qquad{\rho_{F(M)}^{C^{\prime}}}:=\phi_{2}^{-1}(M,C)F({\rho_{M}^{C}}).
\]

\end{claim}

\noindent The proposition bellow is a restatement of \cite[Proposition
4.21]{Ar1}, from which we have kept only the part that we need to prove
Theorem \ref{teo: FS}.

\begin{proposition}
\label{pro coseparability of F gen}Let ${\mathcal{M}}$, $\mathcal{M^{\prime}}%
$, $C$, $C^{\prime}$, $F$ and $F^{\prime}$ be as in
(\ref{ex: co monoidal functor}). We assume that ${\mathcal{M}}$ and
$\mathcal{M^{\prime}}$ are coabelian monoidal categories.

a) If $C$ is coseparable in ${\mathcal{M}}$ then $C^{\prime}$ is coseparable
in $\mathcal{M^{\prime}}$; the converse is true whenever $F^{\prime}$ is separable.

b) Assume that $F$ preserves cokernels. If $C$ is formally smooth
as a coalgebra in ${\mathcal{M}}$ then $C^{\prime}$ is formally
smooth as a coalgebra in ${\mathcal{M^{\prime}}}$; the converse is
true whenever $F^{\prime}$ is separable.
\end{proposition}

\noindent Now we can prove one of the main results of this section.

\begin{theorem}
\label{teo: FS}Let $B$ be a Hopf algebra in a braided abelian and coabelian
monoidal category $(\mathfrak{M},\otimes,\mathbf{1},c)$. We have that:

a) $B$ is coseparable in $(\mathfrak{M}_{B},\otimes,\mathbf{1})$ if and only
if $B$ is coseparable in $(\mathfrak{M},\otimes,\mathbf{1})$.

b) $B$ is formally smooth as a coalgebra in $(\mathfrak{M}_{B},\otimes
,\mathbf{1})$ if and if  $B$ is formally smooth as a coalgebra in
$(\mathfrak{M},\otimes,\mathbf{1})$.
\end{theorem}

\begin{proof}
We apply Proposition \ref{pro coseparability of F gen} in the case when
${\mathcal{M}}:=(\mathfrak{M}_{B},\otimes,\mathbf{1})$ and $\mathcal{M^{\prime
}}:=(\mathfrak{M},\otimes,\mathbf{1})$, which are coabelian categories. We
take $(F,\phi_{0},\phi_{2} ):(\mathfrak{M}_{B},\otimes,\mathbf{1}%
)\rightarrow(\mathfrak{M},\otimes,\mathbf{1})$ to be the forgetful
functor from $\mathfrak{M}_{B}$ to $\mathfrak{M}$, where
$\phi_{0}=\mathrm{Id}_{\mathbf{1}}$ and, for any
$U,V\in\mathfrak{M} _{B}$, we have
$\phi_{2}(U,V)=\mathrm{Id}_{U\otimes V}$. We also take
$F^{\prime}$ to be the forgetful functor from
$^{B}\mathfrak{M}_{B}^{B}$ to $^{B}\mathfrak{M}^{B}$. Since
$(\mathfrak{M} ,\otimes,\mathbf{1})$ is an abelian monoidal
category, then the functor $(-)\otimes
B:\mathfrak{M}\rightarrow\mathfrak{M}$ is additive and right
exact. Hence $F$ preserves cokernels, see \cite[Theorem
3.6]{Ar2}). Thus, in view of Proposition \ref{pro coseparability
of F gen}, to conclude the proof of the theorem, it is enough to
show that $F^{\prime}$ is a separable functor.

For each $\left(  M,{^{B}\!\!\rho_{M}},{\rho_{M}^{B}}\right)  \in
{^{B}{\mathfrak{M}^{B}}}$ we define $\left(  M^{co\left(  B\right)
},i_{M^{co\left(  B\right)  }}\right)  $ to be the equalizer of the maps%
\[
\rho_{M}^{B}:M\rightarrow M\otimes B\qquad\text{and}\qquad\left(  M\otimes
u_{B}\right)  r_{M}^{-1}:M\rightarrow M\otimes B.
\]
Since $B$ is right flat ($\mathfrak{M}$ is an abelian monoidal category), we
can apply the dual version of \cite[Proposition 3.3]{Ar2} to show that
$\left(  M^{co\left(  B\right)  },i_{M^{co\left(  B\right)  }}\right)  $
inherits from $M$ a natural left $B$-comodule structure ${^{B}\!\!\rho
_{M^{co\left(  B\right)  }}}:M^{co\left(  B\right)  }\rightarrow B\otimes
M^{co\left(  B\right)  }$. As a matter of fact, with respect to this comodule
structure, $M^{co\left(  B\right)  }$ is the kernel of $\rho_{M}^{B}-\left(
M\otimes u_{B}\right)  r_{M}\ ^{-1}$ in the category ${^{B}{\mathfrak{M}}}$.
We obtain a functor $F^{\prime\prime}:{}^{B}\mathfrak{M}^{B}\to{}%
^{B}\mathfrak{M}$ defined by:
\[
F^{\prime\prime}( M,{^{B}\!\!\rho_{M}},{\rho_{M}^{B}}) =( M^{co( B) }%
,{^{B}\!\!\rho_{M^{co( B) }}}) .
\]
Then $F^{\prime\prime}\circ F^{\prime}$ associates to $(
M,{\mu_{M}^{B}} ,{^{B}\!\!\rho_{M}},{\rho_{M}^{B}})$ the left
$B$-comodule $( M^{co( B) },{^{B}\!\!\rho_{M^{co( B) }}})$ in
$\mathfrak{M}$. By the dual version of \cite[Proposition
3.6.3]{BD}, it results that $F^{\prime\prime}\circ F^{\prime}$ is
a monoidal equivalence. Therefore $F^{\prime\prime}\circ
F^{\prime}$ is a separable functor and hence $F^{\prime}$ is
separable too.
\end{proof}

\noindent A convenient way to check that a Hopf algebra $B$ is
coseparable in $_{B}\mathfrak{M}$ is to show that $B$ has a total
integral in $\mathfrak{M}$. This characterization of coseparable
Hopf algebras will be proved next.

\begin{definition}
Let $B$ be a Hopf algebra in a braided abelian and coabelian
monoidal category $(\mathfrak{M},\otimes,\mathbf{1},c)$. A
morphism $\lambda:B\rightarrow \mathbf{1}$ in $\mathfrak{M}$ is
called a (left)
\emph{total integral }if it satisfies the relations:%
\begin{align}
r_{B}(B\otimes\lambda)\Delta &  =u\lambda,\label{int1}\\
\lambda u  &  =\operatorname{Id}_{\mathbf{1}}. \label{int2}%
\end{align}

\end{definition}

\begin{figure}[h]
\begin{center}
\fbox{\includegraphics{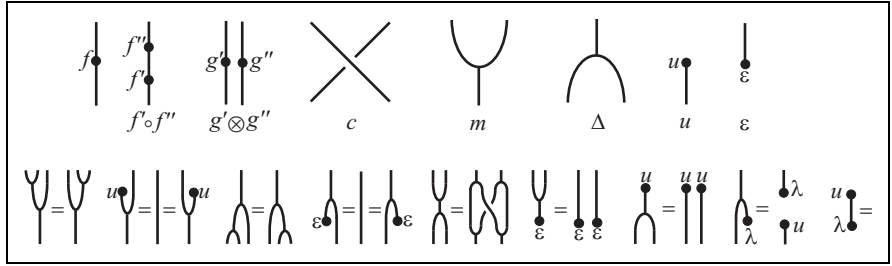}}
\end{center}
\caption{Diagrammatic representation of morphisms in $\mathfrak{M}.$}%
\label{fig:morfisme}%
\end{figure}

\begin{noname}
In order to simplify the computation we will use the diagrammatic
representation of morphisms in a braided category. For details on
this method the reader is referred to \cite[XIV.1]{Ka}. On the
first line of pictures in Figure \ref{fig:morfisme} are included
the basic examples: the representation of a morphism
$f:V\rightarrow W$ (downwards, the domain up) and the diagrams of
$f'\circ f''$ , $g'\otimes g''$ and $c_{V,W}.$ The last four
diagrams denote respectively the multiplication, the
comultiplication, the unit and the counit of a bialgebra $B$ in
$\mathfrak{M}$. The graphical representation of associativity,
existence of unit, coassociativity, existence of counit,
compatibility between multiplication and comultiplication, the
fact that $\varepsilon$ is a morphism of algebras and the fact
that $u$ is a morphism of coalgebras can be found also in
Figure~\ref{fig:morfisme} (second line). The last two pictures on
the same line are equivalent to the definition of a total
integral. The fact that the right hand side of the last equality
is empty means that we can remove the left hand side in any
diagrams that contains it.
\end{noname}

\begin{proposition}
A Hopf algebra $B$ in an abelian and coabelian braided monoidal category
$\mathfrak{M}$ is coseparable in $\mathfrak{M}_{B}$ if and only if it has a
total integral. In this case, $B$ is formally smooth as a coalgebra in
$\mathfrak{M}_{B}$.
\end{proposition}

\begin{figure}[h]
\begin{center}
\fbox{\includegraphics{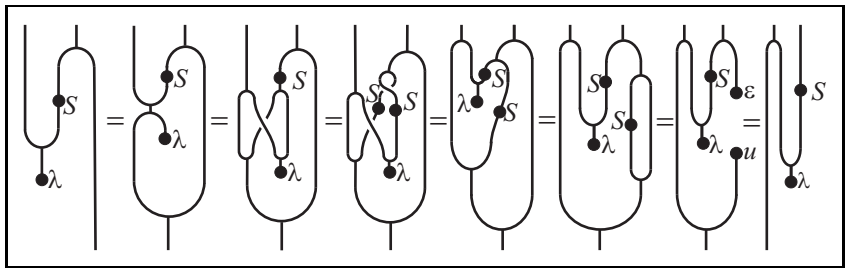}}
\end{center}
\caption{The proof of relation (\ref{ec:lambda}).}%
\label{fig: lambda}%
\end{figure}

\begin{figure}[h]
\begin{center}
\fbox{\includegraphics{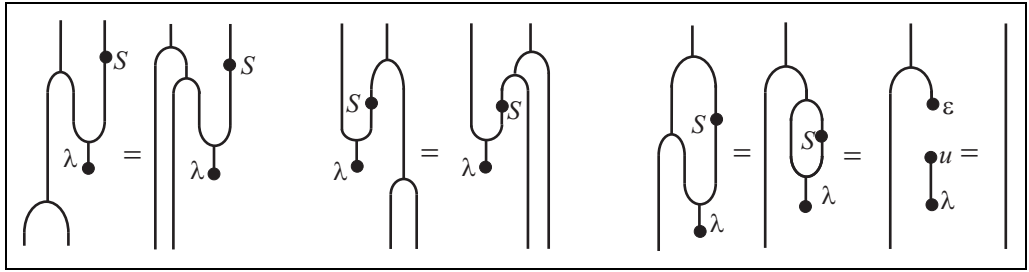}}
\end{center}
\caption{The map $\theta$ is a $B$-bicolinear section of $\Delta.$}%
\label{fig:Bcolin}%
\end{figure}

\begin{figure}[h]
\begin{center}
\fbox{\includegraphics{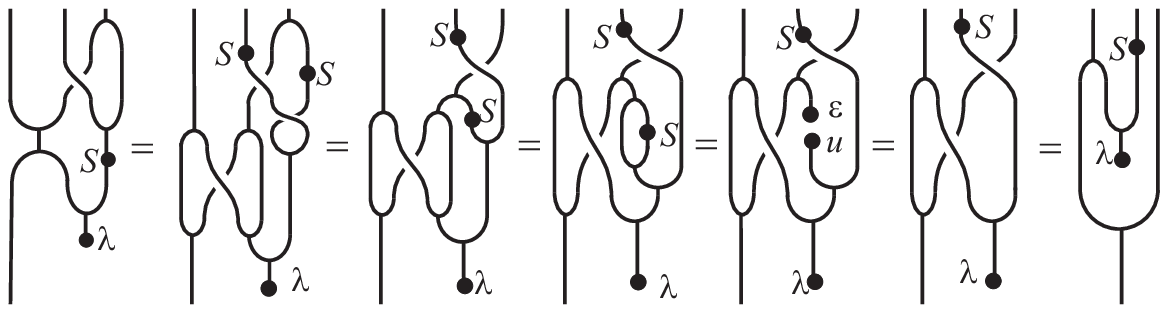}}
\end{center}
\caption{The map $\theta$ is $B$-linear.}%
\label{fig:B_linear}%
\end{figure}

\begin{proof}
We first assume that there is a total integral
$\lambda:B\rightarrow\mathbf{1}$. Let us show that:
\begin{equation}
l_{B}(\lambda m\otimes B)(B\otimes S\otimes B)(B\otimes\Delta)=r_{B}%
(B\otimes\lambda m)(B\otimes B\otimes S)(\Delta\otimes B). \label{ec:lambda}%
\end{equation}
The proof is given in Figure \ref{fig: lambda}. The first equality
follows by relation (\ref{int1}) and the definition of $u$. The
second relation is a consequence of the fact that $B$ is a
bialgebra in $\mathfrak{M},$ so $\Delta m=(m\otimes c_{B,B}\otimes
m)(\Delta\otimes\Delta).$ The third equation follows by the fact
that the antipode $S$ is an anti-morphism of coalgebras, i.e.
$\Delta S=(S\otimes S)c_{B,B}\Delta$. For the fourth equality we
used that the braiding is a functorial morphisms (thus $m,$ $S$
and $\lambda$ can by pulled along the string over and under any
crossing). The last two equalities follow  $m(S\otimes B)\Delta
=u\varepsilon$ and the properties of $\varepsilon$ and $u.$

We now define $\theta:B\otimes B\rightarrow B$ by $\theta(x\otimes
y)=l_{B}(\lambda m\otimes B)(B\otimes S\otimes B)(B\otimes\Delta).$ We have to
prove that $\theta$ is a section of $\Delta$ in the category of $B$%
-bicomodules in $\mathfrak{M}_{B}.$ Note that the category of $B$-bicomodules
in $\mathfrak{M}_{B}$ is $^{B}\mathfrak{M}_{B}^{B}.$ An object $M\in
\mathfrak{M}$ is in $^{B}\mathfrak{M}_{B}^{B}$ if it is a right $B$-module and
a $B$-bicomodule such that $M$ is a Hopf module both in $^{B}\mathfrak{M}_{B}$
and $\mathfrak{M}_{B}^{B}.$

Let us show that $\theta$ is a $B$-bicolinear section of $\Delta.$ Taking into
account relation (\ref{ec:lambda}), we prove that $\theta$ is left
$B$-colinear in the first equality from Figure \ref{fig:Bcolin}. The fact that
$\theta$ is right $B$-colinear is proved in the second equality of the same
figure. In both of them, we used that $\Delta$ is coassociative and that the
comodule structures on $B$ and $B\otimes B$ are defined by $\Delta$,
$B\otimes\Delta$ and $\Delta\otimes B$. To show that $\theta$ is a section of
$\Delta$ we use that $\lambda u=\operatorname{Id}_{\mathbf{1}},$ see the last
sequence of equalities from Figure \ref{fig:Bcolin}.

It remains to prove that $\theta$ is right $B$-linear. This is
done in Figure \ref{fig:B_linear}. The first equality was obtained
by using the fact that $\Delta$ is a morphism of algebras and that
$S$ is an anti-morphism of algebras, i.e. $mS=(S\otimes
S)c_{B,B}m$. To get the second equality we pulled $S$ and $\Delta$
under a crossing (this is possible because the braiding is
functorial). For the third equality we used associativity and
coassociativity. The fourth and the fifth equalities result by the
definitions of the antipode, unit and counit. To deduce the sixth
equality we pulled $m$ and $\lambda$ over the crossing.

Conversely, let $\theta:B\otimes B\rightarrow B$ be a section of
$\Delta,$ which is a morphism of $B$-bicomodules in
$\mathfrak{M}_{B}.$ Let $\lambda:=\varepsilon\theta(B\otimes
u)r_{B}^{-1}.$ Since $\theta$ is a morphism of right $B$-comodules
it results that
\[
\Delta\theta(B\otimes u)r_{B}^{-1}=[\theta(B\otimes u)r_{B}^{-1}\otimes
u]r_{B}^{-1}.
\]
Then, by applying $\varepsilon\otimes B,$ we get $\theta(B\otimes u)r_{B}%
^{-1}=u\lambda.$ As $\theta$ is $B$-colinear, we have:%
\[
\Delta\theta(B\otimes u)r_{B}^{-1}=(B\otimes\theta)(\Delta\otimes u)r_{B}%
^{-1}.
\]
Therefore, by the definition of $\lambda,$ we get (\ref{int1}).
Since $\theta$ is a section of $\Delta$ we deduce that $\lambda
u=\mathrm{Id}_{\mathbf{1}}.$
\end{proof}
\newpage\vspace*{0.1mm}
\vspace*{0.4cm}

\vspace*{0.4cm}

\vspace*{0.4cm}

\vfill\strut\newpage\begin{definition}
Let $E$ be a bialgebra in a braided monoidal category $((\mathfrak{M}%
,\otimes,\mathbf{1}),c)$. Let $H$ be a Hopf subalgebra of $E$. Following
\cite[Definition 5.1]{Schauenburg1}, we say that $\pi:E\rightarrow H$ is a
\emph{(right) weak projection} (onto $H$) if it is a right $H$-linear
coalgebra homomorphism such that $\pi\sigma=\mathrm{Id}_{H},$ where
$\sigma:H\rightarrow E$ is the canonical morphism.
\end{definition}

\noindent The wedge product in an arbitrary braided category was constructed
in \cite[4.7]{AMS1}, while direct systems of Hochschild extensions were
defined in \cite[Definition 4.11]{AMS1}.

\begin{theorem}
\label{teo: Magnum} Let $(\mathfrak{M},\otimes,\mathbf{1},c)$ be a
braided monoidal category. Assume that $\mathfrak{M}$ is
semisimple (i.e. every object is projective) abelian and that both
the functors $X\otimes(-)$ and $(-)\otimes X$ from $\mathfrak{M}$
to $\mathfrak{M}$ are additive for every $X\in \mathfrak{M}$. Let
$A$ be bialgebra in $\mathfrak{M}$ and let $B$ be a sub-bialgebra
of $A$ in $\mathfrak{M}.$ Let $\sigma:B\rightarrow A$ denote the
canonical inclusion. Assume that $A$ is the direct limit (taken in
$\mathfrak{M})$ of $(B^{\wedge^{n}})_{n\in\mathbb{N}},$ the
$B$-adic coalgebra filtration$.$ If $B$ has an antipode (i.e. it
is a Hopf algebra in $\mathfrak{M}$) and $B$ is formally smooth as
a coalgebra in $\mathfrak{M}$ (e.g. $B$ is coseparable in
$\mathfrak{M}$), then $A$ has a right weak projection onto $B$.
\end{theorem}

\begin{proof}
Let $\Delta$ and $\varepsilon$ be respectively the comultiplication and counit
of $A.$ We denote $B^{\wedge^{n+1}}$ by $A_{n},$ where $A_{0}=B.$ We will
denote the canonical projection onto $A/A_{n}$ by $p_{n}.$ One can regard $B$
and $A$ as a coalgebras in $\mathfrak{M}_{B},$ the latter object being a right
$B$-bimodule via $\sigma.$ By induction, it results that $A_{n}\in
\mathfrak{M}_{B},$ for every $n\in\mathbb{N},$ as $A_{n+1}=\operatorname*{Ker}%
[(p_{n}\otimes p_{0})\Delta]$ and by induction hypothesis $p_{n}$ is
$B$-linear (of course $\Delta$ is a morphism of right $B$-modules, since
$\sigma$ is a bialgebra map). Therefore, $(B^{\wedge^{n+1}})_{n\in\mathbb{N}}$
is the $B$-adic filtration on $A$ in $\mathfrak{M}_{B}.$ We want to prove that
the canonical injections $A_{n}\rightarrow A_{n+1}$ split in $\mathfrak{M}%
_{B}$ so that $(A_{n})_{n\in\mathbb{N}}$ is a direct system of
Hochschild extensions in $\mathfrak{M}_{B}$. Indeed, it is enough
to show that the canonical projection $A_{n+1}\rightarrow
A_{n+1}/A_{n}$ has a section in
$\mathfrak{M}_{B}.$ By \cite[Lemma 2.19]{AMS3} it follows that $A_{n+1}%
/A_{n}=(A_{n}\wedge B)/B$ has a canonical right comodule structure $\rho
_{n}:A_{n+1}/A_{n}\rightarrow A_{n+1}/A_{n}\otimes B,$ which is induced by the
comultiplication of $A.$ If $\mu_{n}:A_{n+1}/A_{n}\otimes B\rightarrow
A_{n+1}/A_{n}$ denotes the right $B$-action then $A_{n+1}/A_{n}$ is a
right-right Hopf module, that is:%
\[
\rho_{n}\mu_{n}=(\mu_{n}\otimes m_{B})(M\otimes c_{M,B}\otimes B)(\rho
_{n}\otimes\Delta_{B}).
\]
The structure theorem for Hopf modules works for Hopf algebras in
abelian braided categories, thus
$\mathfrak{M}_{B}^{B}\simeq\mathfrak{M}.$ In fact, the proof of
\cite[Theorem 4.11]{Sw} can be easily written using the graphical
calculus explained above, so it holds in an arbitrary abelian
braided category (note that a similar result, in a braided
category with splitting idempotents can be found in \cite{BGS}).
The equivalence of categories is established by the functor that
associates to $V\in\mathfrak{M}$ the Hopf module $V\otimes B$
(with the coaction $V\otimes\Delta_{B}$ and action $V\otimes
m_{B}$). We deduce that there is a $V\in\mathfrak{M}$ such that
$A_{n+1}/A_{n}\simeq V\otimes B$ (isomorphism in
$\mathfrak{M}_{B}^{B}).$ Thus, to prove that the inclusion
$A_{n}\rightarrow A_{n+1}$ splits in $\mathfrak{M}_{B},$ is
sufficient to show that $V\otimes B$ is projective in
$\mathfrak{M}_{B}$, for any object $V$ in $\mathfrak{M.}$ But, by
\cite[Proposition 1.6]{AMS1}, we
have the adjunction:%
\[
\operatorname{Hom}_{\mathfrak{M}_{B}}(V\otimes B,X)\simeq\operatorname{Hom}%
_{\mathfrak{M}}(V,X),\ \ \ \forall X\in\mathfrak{M}_{B}.
\]
Since, by assumption, $\mathfrak{M}$ is semisimple, we deduce that
$V\otimes B$
is projective in $\mathfrak{M}_{B},$ and hence that $(A_{n}%
)_{n\in\mathbb{N}}$ is a directed system of Hochschild extensions
in $\mathfrak{M}_{B}.$ As $B$ is formally smooth as a coalgebra in
$\mathfrak{M}$, by Theorem \ref{teo: FS}, it is also formally
smooth as a coalgebra in
$\mathfrak{M}_{B}$. Since $A$ is the direct limit of $(A_{n})_{n\in\mathbb{N}%
},$ we conclude the proof by applying \cite[Theorem 4.16]{AMS1} to the case
$\mathcal{M}=\mathfrak{M}_{B}$.
\end{proof}

\noindent As a consequence of Theorem \ref{teo: Magnum} we recover the
following result.

\begin{corollary}
\cite[Theorem 7.38]{Ar1} \label{coro: fs and linear retraction}%
\label{teo: weak}Let $K$ be a field and let $A$ be a bialgebra in
$\mathcal{M}_{K}$ (ordinary bialgebra). Suppose that $B$ is a
sub-bialgebra of $A$ with antipode. Assume that $B$ is formally
smooth as a coalgebra and that $Corad(A)\subseteq B$. Then $A$ has
a right weak projection onto $B$.
\end{corollary}

\noindent Recall that, for a right $H$-comodule $(M,\rho),$ the subspace of
coinvariant elements $M^{co(H)}$ is defined by setting $\ M^{co(H)}%
=\{m\in M\mid\rho(m)=m\otimes1\}.$ If $A$ is an algebra in $\mathcal{M}^{H}$
then $A^{co(H)}$ is a subalgebra of $A$.

When $H$ is a cosemisimple coquasitriangular Hopf algebra, then
$\mathcal{M}^{H}$ is a semisimple braided monoidal category. Note
that bialgebras in $\mathcal{M}^{H}$ are usual coalgebras, so we
can speak about the coradical of a bialgebra in this category.

\begin{corollary}
Let $H$ be a cosemisimple coquasitriangular Hopf algebra and let
$A$ be a bialgebra in $\mathcal{M}^{H}$. Let $B$ denote the
coradical of $A$. Suppose that $B$ is a sub-bialgebra of $A$ (in
$\mathcal{M}^{H}$) with antipode. If $B\subseteq A^{co(H)}$ then
there is a right weak projection $\pi:A\to B$ in
$\mathcal{M}^{H}$.
\end{corollary}

\begin{proof}
Since $B\subseteq A^{co(H)}$ it follows that $c_{B,B}$ is the
usual flip map and $B$ is an ordinary cosemisimple Hopf algebra. A
Hopf algebra is cosemisimple if and only if there is
$\lambda:B\rightarrow K$ such that (\ref{int1}) and (\ref{int2})
hold true (see e.g. \cite[Exercise 5.5.9]{DNR}). Obviously
$\lambda$ is a morphism of $H$-comodules, as $B\subseteq
A^{co(H)}$. The conclusion follows by Theorem \ref{teo: Magnum}.
\end{proof}

\section{Weak Projections onto a Braided Hopf Algebra}

Our main aim in this section is to characterize bialgebras in a braided
monoidal category with a weak projection onto a Hopf subalgebra.

\begin{noname}
\label{assumptions} Throughout this section we will keep the
following assumptions and notations.

\begin{enumerate}
\item[1)] $(\mathfrak{M},\otimes,\mathbf{1},c)$ is an abelian and
coabelian braided monoidal category;

\item[2)] $\left(  A,m_{A},u_{A},\Delta_{A},\varepsilon_{A}\right)
$ is a bialgebra in $\mathfrak{M}$;

\item[3)] $\left(  B,m_{B},u_{B},\Delta_{B},\varepsilon_{B}\right)
$ is a sub-bialgebra of $A$ that has an antipode $S_{B}$ (in
particular $B$ is a Hopf algebra in $\mathfrak{M}$);

\item[4)] $\sigma:B\rightarrow A$ denotes the canonical inclusion
(of course, $\sigma$ is a bialgebra morphism);

\item[5)] $\pi:A\rightarrow B$ is a right weak projection onto $B$
(thus $\pi$ is a morphism of coalgebras in $\mathfrak{M}_{B}$,
where $A$ is a right $B$-module via $\sigma$, and
$\pi\sigma=\mathrm{Id}_{B}$);

\item[6)] We define the following three endomorphisms (in $\mathfrak{M}$) of $A$:%
\[
\Phi:=\sigma
S_{B}\pi,\qquad\qquad\Pi_{1}:=\sigma\pi,\qquad\qquad\Pi_{2}:=m_{A}(
A\otimes\Phi) \Delta_{A}.
\]

\end{enumerate}
\end{noname}

\noindent Our characterization of $A$ as a generalized crossed product is
based on the work of Bespalov and Drabant \cite{BD-Cross2}. We start by
proving certain properties of the operators $\Phi$, $\Pi_{1}$ and $\Pi_{2}$.
They will be used later on to show that the conditions in \cite[Proposition
4.6]{BD-Cross2} hold true.

\begin{lemma}
Under the assumptions and notations in (\ref{assumptions}),
$\Pi_{1}$ is a coalgebra homomorphism such that
$\Pi_{1}\Pi_{1}=\Pi_{1}$ and
\begin{equation}
m_{A}\left(  \Pi_{1}\otimes\Pi_{1}\right)  =\Pi_{1} m_{A}\left(  \Pi
_{1}\otimes\Pi_{1}\right) .
\end{equation}

\end{lemma}

\begin{proof}
Obviously $\Pi_{1}$ is a coalgebra homomorphism as $\sigma$ and $\pi$ are so.
Trivially $\Pi_{1}$ is an idempotent, as $\pi\sigma=\mathrm{Id}_{B}$.
Furthermore, we have%
\[
m_{A}\left(  \Pi_{1}\otimes\Pi_{1}\right)  =m_{A}\left(  \sigma\pi
\otimes\sigma\pi\right)  =\sigma m_{B}\left(  \pi\otimes\pi\right)  =\sigma
\pi\left[  \sigma m_{B}\left(  \pi\otimes\pi\right)  \right] = \Pi_{1}
m_{A}\left(  \Pi_{1}\otimes\Pi_{1}\right) ,
\]
so the lemma is proved.
\end{proof}

\begin{figure}[h]
\begin{center}
\fbox{\includegraphics{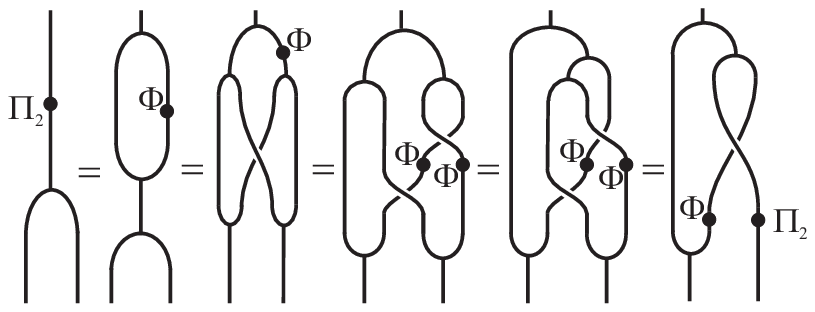} }
\end{center}
\caption{The proof of Eq. (\ref{form: BD1}).}%
\label{fig:BD1}%
\end{figure}

\begin{lemma}
Under the assumptions and notations in (\ref{assumptions}) we
have:
\begin{equation}
\Delta_{A}\circ\Pi_{2}=\left(  m_{A}\otimes A\right)  \circ\left(
A\otimes\Phi\otimes\Pi_{2}\right)  \circ\left(  A\otimes
c_{A,A}\circ\Delta_{A}\right)  \circ\Delta_{A} \label{form: BD1}%
\end{equation}
\end{lemma}
\begin{proof}
See Figure~\ref{fig:BD1} on page \pageref{fig:BD1}. The first equality is
directly obtained from the definition of $\Pi_{2}$. The second equation
follows by the compatibility relation between the multiplication and the
comultiplication of a bialgebra in a braided monoidal category. For the third
equality we used the definition of $\Phi:=\sigma S_{B}\pi$, that $\pi$ and
$\sigma$ are coalgebra homomorphisms and the fact that $S_{B}$ is an
anti-homomorphism of coalgebras in $\mathfrak{M}$, i.e. $\Delta S_{B}%
=(S_{B}\otimes S_{B})c\Delta$. The fourth relation resulted by
coassociativity, while the last one was deduced (in view of
naturality of the braiding) by dragging down one of the
comultiplication morphisms over the crossing and by applying the
definition of $\Pi_{2}$.
\end{proof}

\begin{lemma}
Under the assumptions and notations in (\ref{assumptions}) we
have:
\begin{equation}
\pi\circ\Pi_{2}=u_{B}\circ\varepsilon_{A} \label{form: BD2}%
\end{equation}
\end{lemma}

\begin{figure}[t]
\begin{center}
\fbox{\includegraphics{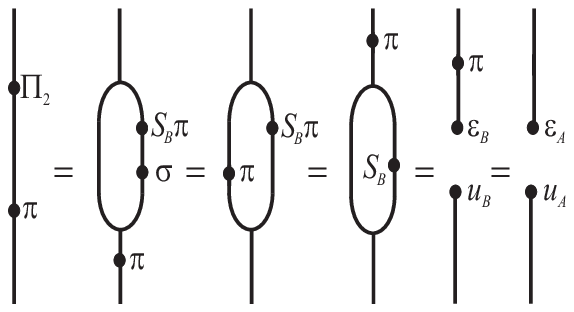}}
\end{center}
\caption{The proof of Eq. (\ref{form: BD2}).}%
\label{fig:BD2}%
\end{figure}

\begin{proof}
See Figure \ref{fig:BD2} on page~\pageref{fig:BD2}. By the
definition of $\Pi_{2}$ we have the first relation. The second one
follows by the fact $\pi$ is right $B$-linear (recall that the
action of $B$ on $A$ is defined by $\sigma$). To deduce the third
equality we use that $\pi$ is a morphism of coalgebras, while the
last relations follow immediately by the properties of the
antipode, unit, counit and $\pi$.
\end{proof}

\begin{figure}[h]
\begin{center}
\fbox{\includegraphics{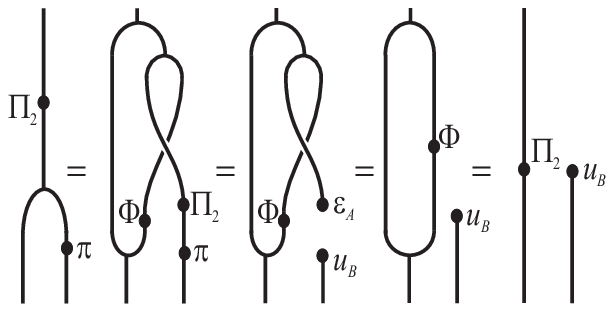}}
\end{center}
\caption{The proof of Eq. (\ref{form: BD5}).}%
\label{fig:BD5}%
\end{figure}

\begin{lemma}
Under the assumptions and notations in (\ref{assumptions}) we
have:
\begin{equation}
\left(  A\otimes\pi\right)  \circ\Delta_{A}\circ\Pi_{2}=\left(  \Pi_{2}\otimes
u_{B}\right)  \circ r_{A}^{-1} \label{form: BD5}%
\end{equation}
\end{lemma}
\begin{proof}
See Figure~\ref{fig:BD5} on page \pageref{fig:BD5}. For the first
equality we used relation (\ref{form: BD1}). The second one
follows by (\ref{form: BD2}), while the third one is a consequence
of the compatibility relation between the counit and the braiding
and the compatibility relation between the counit and the
comultiplication. The last equation is just the definition of
$\Pi_{2}$.
\end{proof}

\begin{lemma}
Under the assumptions and notations in (\ref{assumptions}) we
have:
\begin{equation}
\Pi_{2}m_{A}\left(  A\otimes\sigma\right)  =\Pi_{2} r_{A} \left(
A\otimes\varepsilon_{B}\right)  =r_{A}\left(  \Pi_{2}
\otimes\varepsilon
_{B}\right)  \label{form: BD4}%
\end{equation}

\end{lemma}

\begin{figure}[h]
\begin{center}
\fbox{\includegraphics{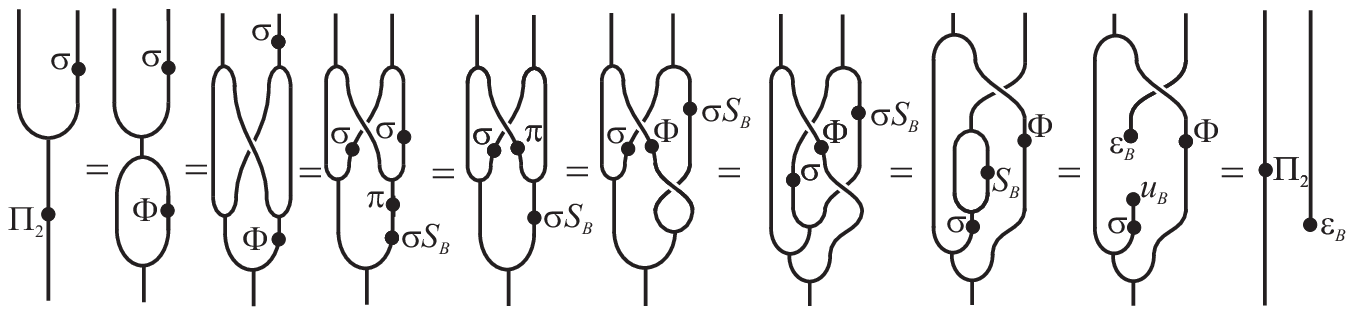}}
\end{center}
\caption{The proof of Eq. (\ref{form: BD4}).}%
\label{fig:BD4}%
\end{figure}

\begin{proof}
The proof can be found in Figure~\ref{fig:BD4} on
page~\pageref{fig:BD4}. The first and the second equalities are
implied by the definition of $\Pi_{2}$ and, respectively, the
compatibility relation between multiplication and comultiplication
in a bialgebra. For the next relation one uses the definition of
$\Phi$ and that $\sigma$ is a morphism of coalgebras. The fourth
relation holds as $\pi$ is right $B$-linear (the $B$-action on $A$
is defined via $\sigma$). The fifth and the sixth equalities
follow as $\sigma$ is a morphism of coalgebras and $S_{B}$ is an
anti-morphism of coalgebras and, respectively, by associativity in
$A$. As the braiding is functorial (so $\Delta _{B}(B\otimes\sigma
S_{B})$ can be dragged under the braiding) and $\sigma$ is a
morphism of algebras we get the seventh equality. To get the
eighth relation we used the definition of the antipode, while the
last one is implied by the properties of the unit and counit in a
Hopf algebra, and the definition of $\Pi_{2}$.
\end{proof}

\begin{lemma}
Under the assumptions and notations in (\ref{assumptions}),
$\Pi_{2}$ is an idempotent such that:
\begin{equation}
\left(  \Pi_{2}\otimes\Pi_{2}\right) \Delta_{A}=\left(  \Pi_{2}\otimes\Pi
_{2}\right) \Delta_{A}\Pi_{2} .\label{form: BD6}%
\end{equation}

\end{lemma}

\begin{figure}[h]
\begin{center}
\fbox{\includegraphics{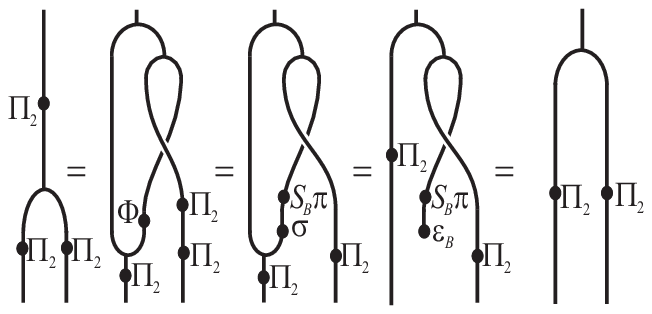}}
\end{center}
\caption{The proof of Eq. (\ref{form: BD6}).}%
\label{fig:BD6}%
\end{figure}

\begin{proof}
Let us first prove that
\begin{equation}
\Pi_{2}\Pi_{2}=\Pi_{2}, \label{form: BD3}%
\end{equation}
that is $\Pi_{2}$ is an idempotent. We have%
\[
\Pi_{2}\Pi_{2} =m_{A}\left(  A\otimes\sigma S_{B}\pi\right)
\Delta_{A} \Pi_{2} =m_{A}\left(  A\otimes\sigma S_{B}\right)
\left(  A\otimes\pi\right) \Delta_{A}\Pi_{2}
\overset{\text{(\ref{form: BD5})}}{=}m_{A}\left( A\otimes\sigma
S_{B}\right)  \left(  \Pi_{2}\otimes u_{B}\right)  r_{A}^{-1}.
\]
Since $\sigma$ and $S_{H}$ are unital morphism and the right unity constraint
$r$ is functorial we get
\[
\Pi_{2}\Pi_{2}=m_{A}\left(  A\otimes u_{A}\right)  \left(  \Pi_{2}%
\otimes\mathbf{1}\right)  r_{A}^{-1} =m_{A}\left(  A\otimes u_{A}\right)
r_{A}^{-1}\Pi_{2}=\Pi_{2}.
\]
For the proof of equation (\ref{form: BD6}) see
Figure~\ref{fig:BD6} on page~\pageref{fig:BD6}. In that figure, we
get the first equality by using (\ref{form: BD1}). The next two
relations are consequences of the definition of $\Phi$, relation
(\ref{form: BD3}) and relation (\ref{form: BD4}). Finally, to
obtain the last equality we use the properties of the antipode and
of the counit, together with the fact that
$\varepsilon_{B}\pi=\varepsilon_{A}$.
\end{proof}

\begin{lemma}
Under the assumptions and notations in (\ref{assumptions}) we
have:
\begin{equation}
\Pi_{1} u_{A}=u_{A} \qquad\text{\emph{and}}\qquad\varepsilon_{A} \Pi
_{2}=\varepsilon_{A}.
\end{equation}

\end{lemma}

\begin{proof}
The first relation is easy: $\Pi_{1} u_{A} =\sigma\pi u_{A}=\sigma\pi\sigma
u_{B}=\sigma u_{B}=u_{A}$. To prove the second one we perform the following
computation:
\[
\varepsilon_{A}\Pi_{2} =\varepsilon_{A} m_{A}\left(  A\otimes\sigma S_{B}%
\pi\right)  \Delta_{A} =m_{\mathbf{1}}\left(  \varepsilon_{A}\otimes
\varepsilon_{A}\right)  \left(  A\otimes\sigma S_{B}\pi\right)  \Delta_{A}.
\]
Therefore, by the properties of the counit of a Hopf algebra in braided
monoidal category, we get
\[
\varepsilon_{A}\Pi_{2} =r_{\mathbf{1}}\left(  \varepsilon_{A}\otimes
\mathbf{1}\right)  \left(  A\otimes\varepsilon_{A}\right)  \Delta
_{A}=\varepsilon_{A} r_{A}\left(  A\otimes\varepsilon_{A}\right)  \Delta
_{A}=\varepsilon_{A},
\]
so the lemma is completely proved.
\end{proof}

\begin{lemma}
Under the assumptions and notations in (\ref{assumptions}), the
homomorphisms
$m_{A}\circ\left(  \Pi_{2}\otimes\Pi_{1}\right)  $ and $\left(  \Pi_{2}%
\otimes\Pi_{1}\right)  \circ\Delta_{A}$ split the idempotent $\Pi_{2}%
\otimes\Pi_{1}$.
\end{lemma}

\begin{figure}[h]
\begin{center}
\fbox{\includegraphics{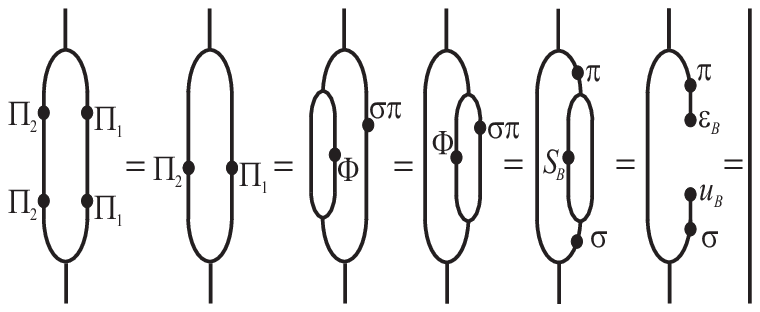}}
\end{center}
\caption{The proof of Eq. (\ref{form: BD12}).}%
\label{fig:BD12}%
\end{figure}

\begin{figure}[h]
\begin{center}
\fbox{\includegraphics{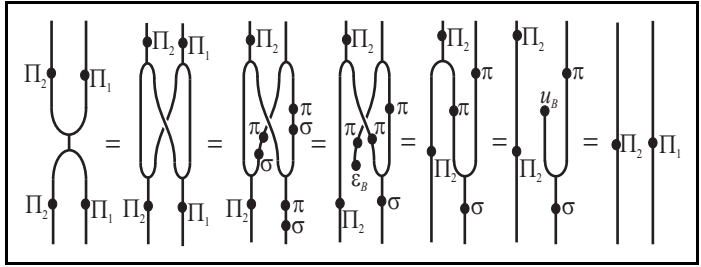}}
\end{center}
\caption{The proof of Eq. (\ref{form: BD13}).}%
\label{fig:BD13}%
\end{figure}

\begin{proof}
We have to prove that
\begin{align}
& \ m_{A}{}\left(  \Pi_{2}\otimes\Pi_{1}\right)  {}\left(  \Pi_{2}\otimes
\Pi_{1}\right)  {}\Delta_{A}=\mathrm{Id}_{A} ,\label{form: BD12}\\
& \left( \Pi_{2}\otimes\Pi_{1}\right)  {}\Delta_{A}{} m_{A}{}\left(  \Pi
_{2}\otimes\Pi_{1}\right)  =\Pi_{1}\otimes\Pi_{2}.\label{form: BD13}%
\end{align}
The proof of relation (\ref{form: BD12}) is shown in
Figure~\ref{fig:BD12} on page~\pageref{fig:BD12}. The first three
equalities are simple consequences of the fact $\Pi_{1}$ and
$\Pi_{2}$ are idempotents, of the definitions of these
homomorphisms and of (co)associativity in $A$. For the fourth
relation we used the definition of $\Phi$ and that $\pi$ is a
morphism of coalgebras and $\sigma$ is a morphism of algebras. The
last two relations result by the definitions of the antipode, unit
and counit in a Hopf algebra, together with
$\varepsilon_{B}\pi=\varepsilon_{A}$ and $\sigma u_{B}=u_{A}$.

The proof of relation (\ref{form: BD13}) is shown in
Figure~\ref{fig:BD13} on page~\pageref{fig:BD13}. The first two
equalities immediately follow by the compatibility relation
between multiplication and comultiplication on $A$, the fact that
$\Pi_{1}$ is a coalgebra homomorphism and $\Pi_{1}=\sigma\pi$. By
using (\ref{form: BD4}) and $\pi
m_A(A\otimes\sigma)=m_B(\pi\otimes B)$, that is $\pi$ is right
$B$-linear, we get the third relation. The relation
$\varepsilon_{B}\pi=\varepsilon_{A}$, the fact that the braiding
is functorial and the properties of the counit are used to obtain
the fourth equality. The fifth one is implied by (\ref{form:
BD5}), while to prove the last equality one uses
$\Pi_{2}\Pi_{2}=\Pi_{2}$, the definition of the unit and the
definition of $\Pi_{1}$.
\end{proof}

\begin{lemma}\label{lem: coinv}
Under the assumptions and notations in (\ref{assumptions}), let
\[
\left(  R,i\right)  :=\mathrm{Eq}\left[  \left( A\otimes\pi\right)  \Delta
_{A},\left(  A\otimes u_{B}\right)  r_{A}^{-1}\right] .
\]
Then there exists a unique morphism $p:A\rightarrow R$ such that%
\[
i p=\Pi_{2}\qquad\text{and}\qquad p i=\mathrm{Id}_{R}.
\]

\end{lemma}

\begin{proof} First of all, we have
\[
\left(  A\otimes\pi\right)
\circ\Delta_{A}\circ\Pi_{2}\overset{(\ref{form: BD5})}{=}\left(
\Pi_{2}\otimes u_{H}\right)  \circ r_{A}^{-1}=\left(  A\otimes
u_{H}\right) \circ\left( \Pi_{2}\otimes\mathbf{1}\right)  \circ
r_{A}^{-1}=\left(  A\otimes u_{H}\right)  \circ
r_{A}^{-1}\circ\Pi_{2}.
\]
Thus, by the universal property of the equalizer, there is a
unique morphism
$p:A\rightarrow R$ such that $ip=\Pi_{2}$. We have%
\begin{align*}
i{} p{} i  &  =\Pi_{2}{} i =m_{A}{}( A\otimes\sigma S_{B}\pi) {}\Delta_{A}{} i
=m_{A}{}( A\otimes\sigma S_{B}) {}( A\otimes\pi) {}\Delta_{A}{} i\\
&  = m_{A}{}( A\otimes\sigma S_{B}) {}( A\otimes u_{B}) {} r_{A}^{-1}{}
i=m_{A}{}( A\otimes u_{A}) {} r_{A}^{-1}{} i =i.
\end{align*}
Since $i$ is a monomorphism we get $p i=\mathrm{Id}_{R}$ so that
$i$ and $p$ split the idempotent $\Pi_{2}.$
\end{proof}

\newpage{}\strut

\vfill\vfill

\vfill\vspace{0.5cm}

\vfill\vspace{0.5cm}

\vfill\strut

\newpage

\begin{noname}
\label{more notation} Before proving one of the main results of
this paper, Theorem~\ref{te: bialgebras}, we introduce some more
notations and terminology. First of all the object $R$, that we
introduced in Lemma \ref{lem: coinv}, will be called \emph{the
diagram of $A$}. Note that $R$ is the `coinvariant subobject' of
$A$ with respect to the right $B$-coaction induced by the
coalgebra homomorphism $\pi$.

We now associate to the weak projection $\pi$ the following data:
\begin{eqnarray*}
 m_{R}:R\otimes R\rightarrow R ,&  &   m_{R}:=p{} m_{A}{}\left(
i\otimes
i\right) ;\\
u_{R}:\mathbf{1}\rightarrow R,&  &   u_{R}:=p{} u_{A};\\
\Delta_{R}:R\to R\otimes R, &  &   \Delta_{R}=(p\otimes p)\Delta_{A}i;\\
\varepsilon_{R}:R\to\mathbf{1},&   &  \varepsilon_{R}=\varepsilon_{A}i;\\
\xi:R\otimes R\rightarrow B,&   &   \xi:=\pi{} m_{A}{}\left(
i\otimes i\right)
;\\
{^{B}\mu_{R}}:B\otimes R\rightarrow R,&   &   {^{B}\mu_{R}}:=p{}
m_{A}{}\left(
\sigma\otimes i\right) ;\\
\mu_{R}^{B}:R\otimes B\rightarrow R,&   &  \mu_{R}^{B}:=p{}
m_{A}{}\left(
i\otimes\sigma\right) ;\\
{^{B}\!\!\rho_{R}} :R\rightarrow B\otimes R,&   &
{^{B}\!\!\rho_{R}}:=\left(
\pi\otimes p\right)  {}\Delta_{A}{} i;\\
 \mu_{B}^{R}:B\otimes R\rightarrow B, & &   \mu_{B}^{R}:=\pi{}
m_{A}{}\left( \sigma\otimes i\right) .
\end{eqnarray*}
\end{noname}

\begin{theorem}
\label{te: bialgebras} We keep the assumptions and notations in
(\ref{assumptions}) and (\ref{more notation}).

1) The diagram $R$ is a coalgebra with comultiplication $\Delta_{R}$ and
counit $\varepsilon_{R}$, and $p$ is a coalgebra homomorphism.

2) The morphisms $m_{A}( i\otimes\sigma)$ and $( p\otimes\pi)\Delta_{A} $ are
mutual inverses, so that $R\otimes B$ inherits a bialgebra structure which is
the cross product bialgebra $R\rtimes B$ defined by\vspace*{-0.4cm}%

\begin{align*}
\hspace*{8mm} m_{R\rtimes B}\ \ = \ \  &  \ ( R\otimes m_{B})( m_{R}\otimes
\xi\otimes B) ( R\otimes{}^{B}\!\mu_{R}\otimes R\otimes R\otimes m_{B})(
R\otimes B\otimes c_{R,R}\otimes R\otimes B\otimes B)\\
&  \ ( R\otimes{}^{B}\!\rho_{R}\otimes R\otimes R\otimes B\otimes B) {}(
\Delta_{R}\otimes\Delta_{R}\otimes B\otimes B)(R\otimes{}^{B}\!\mu_{R}%
\otimes\mu_{B}^{R}\otimes B)\\
& \ ( R\otimes B\otimes c_{B,R}\otimes R\otimes B)( R\otimes\Delta_{B}%
\otimes\Delta_{R}\otimes B)
\end{align*}
\vspace*{-0.8cm}
\begin{align*}
u_{R\rtimes B}  &  =\left(  u_{R}\otimes u_{B}\right)  {}\Delta_{\mathbf{1}%
},\\
\Delta_{R\rtimes B}  &  =( R\otimes m_{B}\otimes R\otimes B) {}( R\otimes
B\otimes c_{R,B}\otimes B) {}( R\otimes{}^{B}\!\rho_{R}\otimes B\otimes B) {}(
\Delta_{R}\otimes\Delta_{B}) ,\\
\varepsilon_{R\rtimes B}  &  =  m_{\mathbf{1}}( \varepsilon_{R}\otimes
\varepsilon_{B}) ,
\end{align*}

\end{theorem}

\begin{proof}
By the previous lemmata, $\Pi_1$ and $\Pi_2$ fulfill the
requirements of the right hand version of \cite[Proposition
4.6(2)]{BD-Cross2}. Thus (1) and (3) of the same result hold. In
our case it can be checked that:
\[
( B_{1},p_{1},i_{1})= ( B,\sigma,\pi) \qquad\text{ and } \qquad( B_{2}%
,p_{2},i_{2})= ( R,i,p).
\]
The explicit form of $m_{R\rtimes B}$ and $\Delta_{R\rtimes B}$ is
a right hand version of the one in the fourth box of diagrams in
\cite[Table 2, page 480]{BD-Cross2}.
\end{proof}

We are now going to investigate a particular case of the above theorem.
Namely, when $A$ is cocommutative and $\xi$ is trivial, we will show that $A$
is the double cross product of a matched pair (see definitions bellow).

\begin{definition}
\label{de: matched pair} Let $\left(  R,m_{R},u_{R},\Delta_{R},\varepsilon
_{R}\right)  $ and $\left(  B,m_{B},u_{B},\Delta_{B},\varepsilon_{B}\right)  $
be bialgebras in a braided abelian and coabelian monoidal category
$(\mathfrak{M},\otimes,\mathbf{1},c).$ Following \cite[Definition 7.2.1, page
298]{Maj-Found}, we say that $\left(  R,B\right)  $ defines a \emph{matched
pair of bialgebras}, if there exist morphisms%
\[
\mathbf{\vartriangleright}:B\otimes R\rightarrow R\qquad\text{and}%
\qquad\mathbf{\vartriangleleft}:B\otimes R\rightarrow B
\]
satisfying the seven conditions below:

\begin{enumerate}
\item $( R,\Delta_{R},\varepsilon_{R} ,{\vartriangleright}) $ is a left
$B$-module coalgebra;

\item $( B,\Delta_{B},\varepsilon_{B} ,{\vartriangleleft}) $ is a right
$R$-module coalgebra;

\item ${\vartriangleleft}( u_{B}\otimes R) =u_{B}{} \varepsilon_{R}{} l_{R};$

\item ${{\vartriangleright}}{}( B\otimes u_{R}) =u_{R}{} \varepsilon_{B}{}
r_{B};$

\item $m_{B}({\vartriangleleft} \otimes\,B)( {B} \otimes{\vartriangleright}
\otimes{\vartriangleleft})( B\otimes\Delta_{B\otimes R}) = {{\vartriangleleft
}}( m_{R}\otimes B) ;$

\item $m_{R}{}( R\otimes{\vartriangleright}) ({{\vartriangleright}}%
\otimes{{\vartriangleleft}}\otimes R) {}( \Delta_{B\otimes R}\otimes R)
={{\vartriangleright}}{}( B\otimes m_{R}) ;$

\item $({\vartriangleleft}\otimes{\vartriangleright})
{}\Delta_{B\otimes R}=c_{R,B}{}(
{\vartriangleright}\otimes{\vartriangleleft}) {}\Delta_{B\otimes
R}$.\end{enumerate} In this case, for sake of shortness, we will
say that $(R,B,\mathbf{\vartriangleright ,\vartriangleleft})$
\emph{is a matched pair of bialgebras in}
$(\mathfrak{M},\otimes,\mathbf{1},c).$
\end{definition}

\begin{noname}\label{bicrossed product} Let
$( R,B,\vartriangleright ,\vartriangleleft)$ be a matched pair. By
\cite[Corollary 2.17]{BD-Cross1}, we get that:
\begin{align*}
m_{R\Join B}  &  =( m_{R}\otimes m_{B}) {}( R\otimes{\vartriangleright}%
\otimes{\vartriangleleft}\otimes B) {}( R\otimes B\otimes c_{B,R}\otimes
R\otimes B) {}( R\otimes\Delta_{B}\otimes\Delta_{R}\otimes B) ,\\
u_{R\Join B}  &  =( u_{R}\otimes u_{B}) {}\Delta_{\mathbf{1}},\\
\Delta_{R\Join B}  &  =( R\otimes c_{R,B}\otimes B) {}( \Delta_{R}%
\otimes\Delta_{B}) ,\\
\varepsilon_{R\Join B}  &  =m_{\mathbf{1}}{}( \varepsilon_{R}\otimes
\varepsilon_{B}) .
\end{align*}
defines a new bialgebra $R\Join B$, that is called \emph{the
double cross product bialgebra}. It can be obtained as a
particular case of the cross product bialgebra $R\rtimes B$ by
setting $^B\!\mu_R=\tl$, $\mu^R_B=\tr$ and taking $\xi$,
$^B\!\rho_R$ and $\mu_B^R$ to be trivial in (\ref{more notation}).
\end{noname}

\begin{theorem}
\label{te: matched pair} Let $\sigma:B\to A$ and $i:R\to A$ be
bialgebra morphisms in a braided monoidal category
$(\mathfrak{M},\otimes,\mathbf{1},c)$ such that
$\Phi:=m_{A}(i\otimes \sigma)$ is an isomorphism in
$\mathfrak{M}$. Let $\Psi:=\Phi^{-1}\Theta$, where
$\Theta:B\otimes R\to A$ is defined by $\Theta:=m_{A}(
\sigma\otimes i)$.\\ Consider the homomorphisms
${\vartriangleright}:B\otimes R\to R$ and
$\vartriangleleft:B\otimes R\to B$ defined by:
\begin{equation}\label{eq:lt and tr}
{\vartriangleright}:=r_{R}(R\otimes\varepsilon_{B})\Psi,\qquad
{\vartriangleleft}:=l_{B}(\varepsilon_{R}\otimes B)\Psi
\end{equation}
 Then $(R,B,\tr,\tl)$ is a matched pair and $A\simeq
R\Join B$.
\end{theorem}

\begin{figure}[h]
\begin{center}
\fbox{\includegraphics{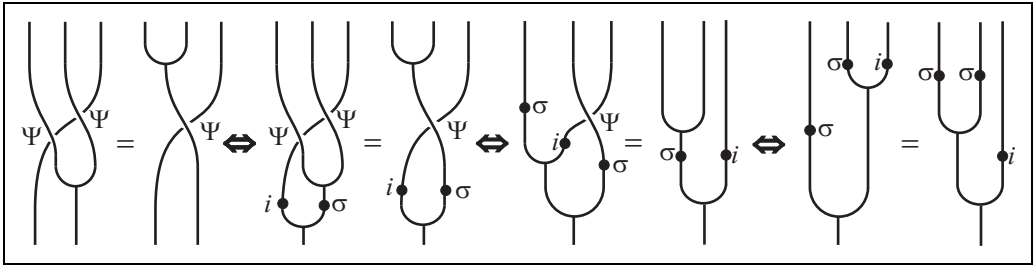}}
\end{center}
\caption{The proof of the first equation in (\ref{form: maj 0'}).}%
\label{fig:cp0}%
\end{figure}

\begin{figure}[h]
\begin{center}
\fbox{\includegraphics{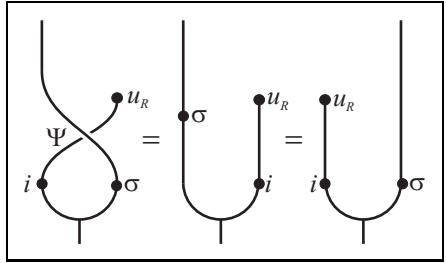}}
\end{center}
\caption{The proof of the second equation in (\ref{form: maj 0'}).}%
\label{fig:cp00}%
\end{figure}

\begin{figure}[h]
\begin{center}
\fbox{\includegraphics{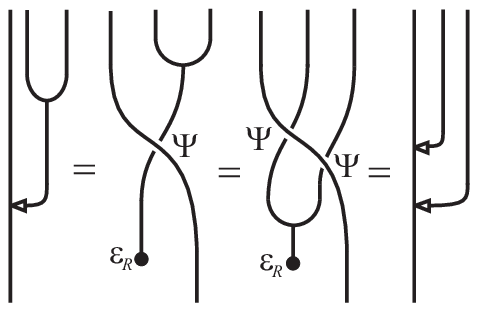}}
\end{center}
\caption{ $B$ is a right $R$-module with respect to $\vartriangleleft$.}%
\label{fig:cpa}%
\end{figure}

\begin{figure}[h]
\begin{center}
\fbox{\includegraphics{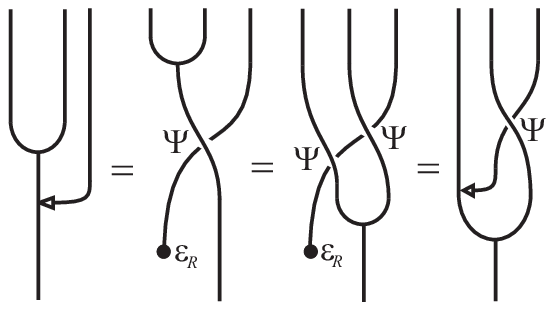}}
\end{center}
\caption{The proof of Eq. (\ref{form: cpb}).}%
\label{fig:cpb}%
\end{figure}

\begin{figure}[h]
\begin{center}
\fbox{\includegraphics{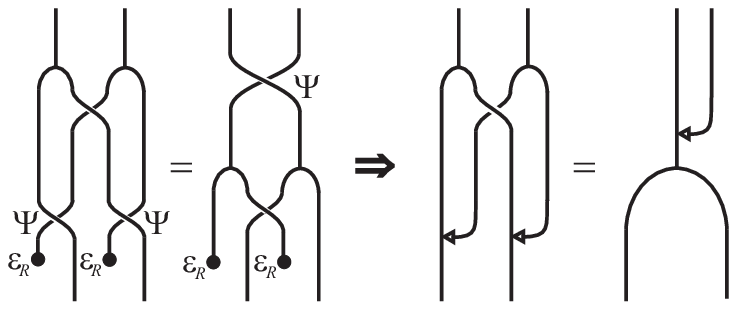}}
\end{center}
\caption{The proof of Eq. (\ref{form: cp1}).}%
\label{fig:cp1}%
\end{figure}

\begin{figure}[h]
\begin{center}
\fbox{\includegraphics{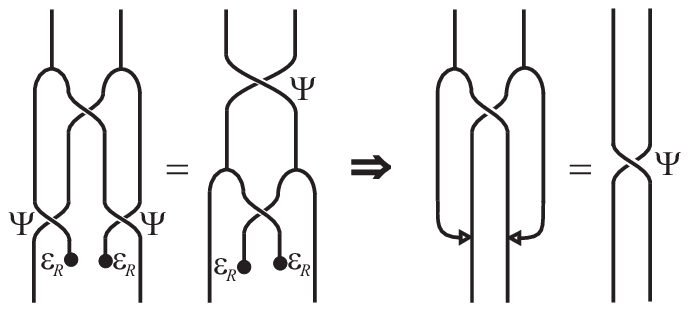}}
\end{center}
\caption{The proof of Eq. (\ref{form: cp2}).}%
\label{fig:cp2}%
\end{figure}

\begin{figure}[hh]
\begin{center}
\fbox{\includegraphics{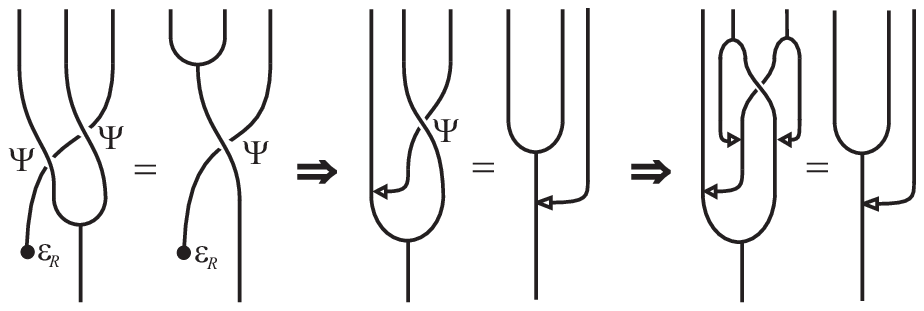}}
\end{center}
\caption{The proof of Eq. (\ref{form: cp3}).}%
\label{fig:cp3}%
\end{figure}

\begin{proof}
We will follow the proof of \cite[Theorem 7.2.3]{Maj-Found}. It is
easy to see that the proofs of relations \cite[(7.10)]{Maj-Found}
and \cite[(7.11)]{Maj-Found} work in a braided monoidal category,
as they can be done in a diagrammatic way. Therefore, we have:
\begin{align}
& ( R\otimes m_{B}) ( \Psi\otimes B) ( B\otimes\Psi) =\Psi( m_{B}\otimes
R),\qquad\Psi(B \otimes u_{R})r_{B}^{-1}=(u_{R}\otimes B)l_{B}^{-1}
.\label{form: maj 0'}\\
& ( m_{R}\otimes B) ( R\otimes\Psi) ( \Psi\otimes R) =\Psi( B\otimes
m_{R}),\qquad\Psi(u_{B}\otimes R)l_{R}^{-1}=(R\otimes u_{B})r_{R}%
^{-1}.\label{form:
maj 0}%
\end{align}
For example the first relation in (\ref{form: maj 0'}) is proved
in Figure~\ref{fig:cp0}. The first equivalence there holds since
$\Phi =m_{A}(i\otimes\sigma)$ is by assumption an isomorphism. The
second and the third equivalences are consequences of
associativity in $A$ and of relation $\Phi\Psi=\Theta$. Since the
last equality is obviously true by associativity in $A$, the
required relation is proved. The second relation in (\ref{form:
maj 0'}) follows by the computation performed in
Figure~\ref{fig:cp00}. The first equality holds since
$\Phi\Psi=\Theta$, while the second results by the fact that $i$
and $\sigma$ are homomorphisms of algebras and by the definition
of the unit in an algebra. To get the second relation in
(\ref{form: maj 0'}) we use the fact that $\Phi$ is an
isomorphism.

As in the proof of \cite[Theorem 7.2.3]{Maj-Found}, by applying
$l_{B}( \varepsilon_{R}\otimes B) $ and $r_{R}(R
\otimes\varepsilon_{R}) $ respectively to (\ref{form: maj 0}) and
(\ref{form: maj 0'}), we get that ${\vartriangleright}$ defines a
left action of $B$ on $R$ and ${\vartriangleleft}$ defines a right
action of $R$ on $B$. Indeed, by applying $\varepsilon_R\otimes B$
to the second relation in (\ref{form: maj 0'}) it is easy to see
that $\tl$ is unital. The second axiom that defines a right action
is checked in Figure~\ref{fig:cpa}.

Furthermore, by applying $l_{B}( \varepsilon_{R}\otimes B) $ and $r_{R}(R
\otimes\varepsilon_{R}) $ respectively to (\ref{form: maj 0'}) and
(\ref{form: maj 0}), we get
\begin{align}
& {\vartriangleleft}(m_{B}\otimes R)=m_{B}({\vartriangleleft}\otimes
B)(B\otimes\Psi)\label{form: cpb}\\
& {\vartriangleright}(B\otimes m_{R})=m_{R}(R\otimes{\vartriangleleft}%
)(\Psi\otimes R)
\end{align}
For the proof of (\ref{form: cpb}) see Figure~\ref{fig:cpb}. We
now want to check that $\Theta: B\otimes R\rightarrow A$ is a
coalgebra homomorphism, where the coalgebra structure on $B\otimes
R$ is given by:
\[
\Delta_{B\otimes R} :=( B\otimes c_{B,R}\otimes R) ( \Delta_{B}\otimes
\Delta_{R}),\qquad\varepsilon_{B\otimes R} :=m_{\mathbf{1}}\left(
\varepsilon_{B}\otimes\varepsilon_{R}\right) .
\]
Indeed, we have
\begin{align*}
\Delta_{A} \Theta &  =( m_{A}\otimes m_{A}) ( A\otimes c_{A,A}\otimes A) (
\Delta_{A}\otimes\Delta_{A}) ( \sigma\otimes i)\\
&  =( m_{A}\otimes m_{A}) ( A\otimes c_{A,A}\otimes A) ( \sigma\otimes
\sigma\otimes i\otimes i) ( \Delta_{B}\otimes\Delta_{R})\\
&  =( m_{A}\otimes m_{A}) ( \sigma\otimes i\otimes\sigma\otimes i) ( B\otimes
c_{B,R}\otimes R) ( \Delta_{B}\otimes\Delta_{R})\\
&  =(\Theta \otimes\Theta)( B\otimes c_{B,R}\otimes R) (
\Delta_{B}\otimes\Delta_{R})=(\Theta
\otimes\Theta)\Delta_{B\otimes R},
\end{align*}
and $\varepsilon_{A} \Theta=m_{\mathbf{1} }( \varepsilon_{A}\otimes
\varepsilon_{A}) ( \sigma\otimes i) =m_{\mathbf{1}}( \varepsilon_{B}
\otimes\varepsilon_{R})=\varepsilon_{B\otimes R}.$ In a similar way, by
interchanging $B$ and $R$, we can prove that $\Phi$ is a homomorphism of
coalgebras. Thus $\Psi=\Phi^{-1}\Theta$ is a coalgebra homomorphism too, so
\begin{equation}
\Delta_{R\otimes B}\Psi=\left(  \Psi\otimes\Psi\right)  \Delta_{B\otimes
R}\qquad\text{and}\qquad\varepsilon_{R\otimes B} \Psi=\varepsilon_{B\otimes
R}. \label{form: maj 1}%
\end{equation}
By applying $l_{B}\left(  \varepsilon_{R}\otimes B\right)  \otimes
l_{B}\left(  \varepsilon_{R}\otimes B\right)  $ to both sides of
the first equality in (\ref{form: maj 1}) we get the first
relation in Figure~\ref{fig:cp1}. By the properties of
$\varepsilon_{B}$ and $\varepsilon_{R}$ we get the second relation
in the same figure, that is we have:
\begin{equation}
\Delta_{B}{\vartriangleleft}=({\vartriangleleft}\otimes{\vartriangleleft})
\Delta_{B\otimes R}.\label{form:
cp1}%
\end{equation}
As $\varepsilon_{B}{\vartriangleleft}=\varepsilon_{B}l_{B}(\varepsilon
_{R}\otimes B)\Psi=\varepsilon_{R\otimes B}\Psi=\varepsilon_{B\otimes R}$ we
have proved that $( B,\Delta_{B},\varepsilon_{B} ,{\vartriangleleft}) $ is a
right $R$-module coalgebra. Similarly one can prove that $( R,\Delta
_{R},\varepsilon_{R} ,{\vartriangleright})$ is a left $B$-module coalgebra.

\begin{figure}[hhh]
\begin{center}
\fbox{\includegraphics{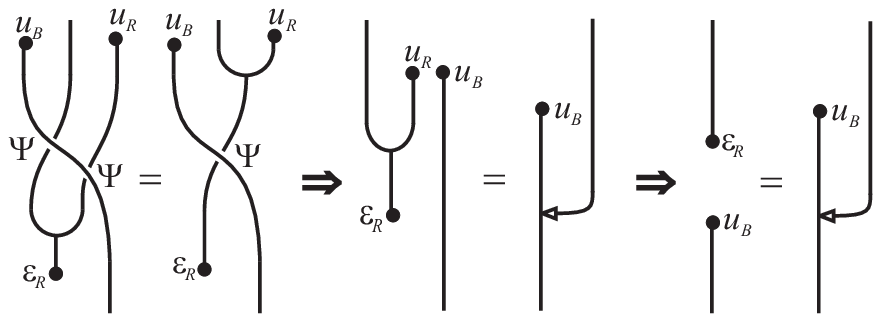}}
\end{center}
\caption{The proof of Eq. (\ref{form: cp4}).}%
\label{fig:cp4}%
\end{figure}

By applying $r_{R}\left(  R\otimes\varepsilon_{H}\right)  \otimes
l_{H}\left( \varepsilon_{R}\otimes H\right)  $ and $l_{H}\left(
\varepsilon_{R}\otimes H\right)  \otimes r_{R}\left(
R\otimes\varepsilon_{H}\right)  $ respectively to both sides of
the first equality in (\ref{form: maj 1}) (see e.g.
Figure~\ref{fig:cp2}) one can prove the relations:
\begin{equation}
\Psi=( {\vartriangleright}\otimes{\vartriangleleft})\Delta_{B\otimes R}%
\qquad\text{and}\qquad c_{R,B}\Psi=( {\vartriangleleft}\otimes
{\vartriangleright}) \Delta_{B\otimes R}. \label{form: cp2}%
\end{equation}
By the two relations of (\ref{form: cp2}) we deduce%
\begin{equation}
c_{R,B}( {\vartriangleright}\otimes{\vartriangleleft})\Delta_{B\otimes R}=(
{\vartriangleleft}\otimes{\vartriangleright}) \Delta_{B\otimes R}.
\label{form: maj Match 7.9}%
\end{equation}
By applying $\varepsilon_{R}\otimes B$ to the both sides of
(\ref{form: maj 0'}) we get the first equation in Figure~\ref{fig:cp3}. By the
definition of the right action of $R$ on $B$ we get the relation in the middle
of that figure. By using the first equation in (\ref{form: cp2}) we get the
last equality in Figure~\ref{fig:cp3}. Therefore we have proved the following
equation:%
\begin{equation}
m_{B}({\vartriangleleft}\otimes B) ( B\otimes{\vartriangleright}%
\otimes{\vartriangleleft}) ( B\otimes\Delta_{B\otimes R}) ={\vartriangleleft}(
m_{R}\otimes B). \label{form: cp3}%
\end{equation}
The relation 6) from Definition~\ref{de: matched pair} can be
proved similarly. Finally, by composing both sides of (\ref{form:
maj 0}) by $\varepsilon_{R}\otimes B$ to the left and by
$u_{B}\otimes B\otimes u_{R}$ to the right we get:
\begin{equation}
{\vartriangleleft}( u_{B}\otimes R) =u_{B}\varepsilon_{R}l_{R}
\label{form: cp4}%
\end{equation}
The details of the proof are given in Figure~\ref{fig:cp4}.
Analogously one can prove relation 4) from Definition~\ref{de:
matched pair}. In conclusion, we have proved that
$(R,B,\vartriangleleft,\vartriangleright)$ is a matched pair and
that $\Phi$ is a morphism of coalgebras.

It remains to prove that $\Phi$ is an isomorphism of algebras.
Obviously $\Phi$ is an unital homomorphism. By (\ref{form: cp2})
it follows that $m_{R\Join B}=(m_{R}\otimes
m_{B})(R\otimes\Psi\otimes B)$. Since $i$ and $\sigma$ are
morphisms of algebras and $m_{A}$ is associative we get
\begin{align*}
\Phi m_{R\Join B}  &  = m_{A}(i\otimes\sigma)(m_{R}\otimes
m_{B})(R\otimes \Psi\otimes B)= m_{A}(m_{A}\otimes A)[i\otimes
m_{A}(i\otimes\sigma
)\Psi\otimes\sigma]\\
&  = m_{A}(m_{A}\otimes A)(i\otimes\Phi\Psi\otimes\sigma) =
m_{A}(m_{A}\otimes
A)(i\otimes\Theta\otimes\sigma)\\
&  = m_{A}(m_{A}\otimes A)[i\otimes m_{A}(\sigma\otimes
i)\otimes\sigma] =m_{A}(\Phi\otimes\Phi).
\end{align*}
Trivially $\Phi u_{R\Join B}=u_{A}$ since $i$ and $\sigma$ are
unital homomorphism and $m_{A}(u_{A}\otimes
u_{A})\Delta_{\mathbf{1}}=u_{A}$, so the theorem is proved.
\end{proof}

\vspace*{0.7cm}

\vspace*{0.7cm}

\vspace*{0.7cm}

\vfill\strut\newpage

\vspace*{0.7cm}

\vspace*{0.7cm}

\vspace*{0.7cm}

\vfill\strut\newpage

\begin{theorem}
\label{te: double cross product} We keep the assumptions and
notations in (\ref{assumptions}) and (\ref{more notation}). We
also assume that $A$ is cocommutative and $\xi$ is trivial, i.e.
$c_{A,A}\Delta_{A}=\Delta_{A}$ and $\xi=u_{B} m_{\mathbf{1}}(
\varepsilon_{R}\otimes\varepsilon_{R})$. Then
\[
\left(  R,B,\vartriangleright,\vartriangleleft\right)
\]
is a matched pair of bialgebras such that $A\simeq R\Join B$,
where
\begin{equation}\label{eq: actions l}
\mathbf{\vartriangleright:=}{^{B}\mu_{R}}=p m_{A}\left(
\sigma\otimes i\right)  :B\otimes R\rightarrow R\qquad\text{and}%
\qquad\mathbf{\vartriangleleft:=}\mu_{B}^{R}=\pi m_{A}\left(
\sigma\otimes i\right)  :B\otimes R\rightarrow B.
\end{equation}
\end{theorem}

\begin{proof}
Since $\xi$ is trivial, by \cite[Proposition 3.7(5)]{BD-Cross2} it
follows that $R$ is an algebra and $i:R\to A$ is an algebra
homomorphism. Our aim now is to show that $i:R\rightarrow A$ is a
coalgebra homomorphism too. In view of \cite[Proposition 3.7
(8)]{BD-Cross2} it is enough to prove that
\begin{equation}
\label{eq:rho}{^{B}\!\rho_{R}}=( u_{B}\otimes R)l_{R}^{-1}\quad
\text{and}\quad (\pi\ot \pi)\Delta_A i=(u_B\ot
u_B)\Delta_{\mathbf{1}}\varepsilon_R.
\end{equation}
Since $\pi$ is a coalgebra homomorphism, the second equality
follows by \cite[Proposition 3.7 (6)]{BD-Cross2}. Let us prove the
first one. Indeed, as $R$ is the equalizer of
$(A\otimes\pi)\Delta_{A}$ and $(A\otimes
u_{B})r_{A}^{-1}$, we have%
\begin{align*}
( p\otimes\pi) \Delta_{A} i  &  =( p\otimes B) ( A\otimes\pi) \Delta_{A} i=(
p\otimes B) ( A\otimes u_{B}) r_{A}^{-1} i =( R\otimes u_{B}) ( p\otimes
\mathbf{1}) r_{A}^{-1} i\\
& =( R\otimes u_{B}) r_{R}^{-1} p i =( R\otimes u_{B}) r_{R}^{-1}.
\end{align*}
Therefore%
\begin{align*}
{^{B}\rho_{R}}  &  =( \pi\otimes p) \Delta_{A} i =( \pi\otimes p)
c_{A,A}\Delta_{A} i =c_{A,A}( p\otimes\pi) \Delta_{A} i =c_{A,A}( R\otimes
u_{B}) r_{R}^{-1}\\
&  =( u_{B}\otimes R) c_{R,\mathbf{1}} r_{R}^{-1} =( u_{B}\otimes R)
l_{R}^{-1}.
\end{align*}
Hence (\ref{eq:rho}) is proved and, in consequence, it results
that $i$ is a morphism of bialgebras. By Theorem \ref{te:
bialgebras}(2) the morphisms $\Phi =m_{A}(i\otimes\sigma)$ and $(
p\otimes\pi)\Delta_{A} $ are mutual inverses. Thus we can apply
Theorem~\ref{te: matched pair}. In our case
\[\Psi=\Phi^{-1}\Theta=( p\otimes\pi)\Delta_{A}m_{A}(\sigma\otimes i).\]
In view of (\ref{eq:lt and tr}), it results
\begin{equation*}{\vartriangleright}=r_R(R\ot\varepsilon_B)(p\ot\pi)\Delta_A
m_A(\sigma\ot i)=pr_A(A\ot\varepsilon_B)\Delta_A m_A(\sigma\ot
i)=pm_A(\sigma\ot i).
\end{equation*}
In a similar way we get
\begin{equation*}
\tr=\pi m_A(\sigma\ot i).
\end{equation*}
\end{proof}

\begin{noname}\label{nn:smash}
We keep the assumptions and the notations in (\ref{assumptions})
and (\ref{more notation}). We take $A$ to be a cocommutative
bialgebra in $\mathfrak{M}$ with trivial cocycle $\xi$. Thus, by
our results, $A$ is the double cross product of a certain matched
pair $(R,B,\vartriangleright,\vartriangleleft)$, where the actions
$\vartriangleright$ and $\vartriangleleft$ are defined by
relations (\ref{eq: actions l}). Our aim now is to investigate
those bialgebras $A$ as above which, in addition, have the
property that the right action $\vartriangleleft:B\otimes R\to B$
is trivial. We will see that in this case the left action
$\vartriangleright :B\otimes R\to R$ is the adjoint action. More
precisely, we have $i\tr=\mathrm{ad}$, where $\mathrm{ad}$ is
defined by:
\[
\mathrm{ad}=m_A(m_A\otimes A)(\sigma\otimes i\otimes \sigma
S_B)(B\otimes c_{B,R})(\Delta_{B}\otimes R).
\]
Moreover, $A$ can be recovered from $R$ and $B$ as the
`bosonization' $R\#B$, that is $A$ is the smash product algebra
between $R$ and $B$, and as a coalgebra $A$ is isomorphic to the
tensor product coalgebra $R\otimes B$. Recall that the
multiplication and the comultiplication on $R\# B$ are given:
\begin{align*}
m_{R\# B}  &  =( m_{R}\otimes m_{B}) {}( R\otimes\tr\otimes
R\otimes B) {}( R\otimes B\otimes c_{B,R}\otimes B) {}(
R\otimes\Delta_{B}\otimes
R\otimes B) ,\\
u_{R\# B}  &  =( u_{R}\otimes u_{B}) {}\Delta_{\mathbf{1}},\\
\Delta_{R\# B}  &  =( R\otimes c_{R,B}\otimes B) {}( \Delta_{R}\otimes
\Delta_{B}) ,\\
\varepsilon_{R\# B}  &  =m_{\mathbf{1}}{}( \varepsilon_{R}\otimes
\varepsilon_{B}) .
\end{align*}
\end{noname}

\begin{proposition}\label{pr: bosonization}
We keep the assumptions and notations in (\ref{assumptions}) and
(\ref{more notation}). We also assume that $A$ is cocommutative
and $\xi$ is trivial.

a) The action $\vartriangleleft:B\otimes R\to B$ is trivial if and only if
$\pi$ is left $B$-linear.

b) If $\vartriangleleft$ is trivial then the left action
$\vartriangleright :B\otimes R\to R$ is the adjoint action.

c) If $\vartriangleleft$ is trivial then $A\simeq R\#B$, where $B$
acts on $R$ by the left adjoint action.
\end{proposition}

\begin{proof}
Since $(R,i)$ is the equalizer of $(A\otimes\pi)\Delta_{A}$ and $(A\otimes
u_{B})r_{A}^{-1}$ we get
\[
(A\otimes\pi)\Delta_{A}i=(A\otimes u_{B})r_{A}^{-1}i
\]
By applying $\varepsilon_{R}\otimes B$ to the both sides of this
relation we get $\pi i=u_{B}\varepsilon_{R}$. Now we can prove a).
If we assume that $\pi$ is left $B$-linear, i.e. $\pi
m_{A}(\sigma\otimes A)=m_{B}(B\otimes\pi )$, then it results
\[
{\vartriangleleft}=\pi m_A(\sigma\otimes i)=m_{B}(B\otimes\pi
i)=m_{B}(B\otimes
u_{B}\varepsilon_{R})=r_{B}(B\otimes\varepsilon_{R}).
\]
This means that ${\vartriangleleft}$ is trivial. Conversely, let us assume
that $\pi m_{A}(\sigma\otimes i)=r_{B}(B\otimes\varepsilon_{R})$. In order to
prove  that $\pi m_{A}(\sigma\otimes A)=m_{B}(B\otimes\pi)$ we compute $\pi
m_{A}(\sigma\otimes\Phi).$ We get%
\begin{align*}
\pi m_{A}(\sigma\otimes\Phi)  & =\pi m_{A}\left[  \sigma\otimes m_{A}%
(\sigma\otimes i)\right]  =\pi m_{A}(\sigma\otimes i)(m_{B}\otimes
R) =r_{B}(B\otimes\varepsilon_{R})(m_{B}\otimes R)\\
&=m_{B}\left[ B\otimes\pi m_{A}(\sigma\otimes i)\right]
=m_{B}(B\otimes\pi\Phi).
\end{align*}
Since $B\otimes\Phi$ is an isomorphism we deduce the required equality.

b) The proof of $i\tr=\mathrm{ad}$ is given in
Figure~\ref{fig:action}. The definition of the action $\tr$
together with $ip=m_A(A\ot\Phi)\Delta_A$ and $\Phi=\sigma S_B\pi$
yield the first equality. The next one is obtained by applying the
compatibility relation between $m_A$ and $\Delta_A$ and the fact
that $\sigma$ is a morphism of coalgebras. By the first part of
the proposition, $\pi$ is left $B$-linear. Thus we have  the third
equality. By using $\pi i=u_B\varepsilon_R$ and the properties of
the unit and counit we conclude the proof of $i\tr=\mathrm{ad}$.

\begin{figure}[hh]
\begin{center}
\fbox{\includegraphics{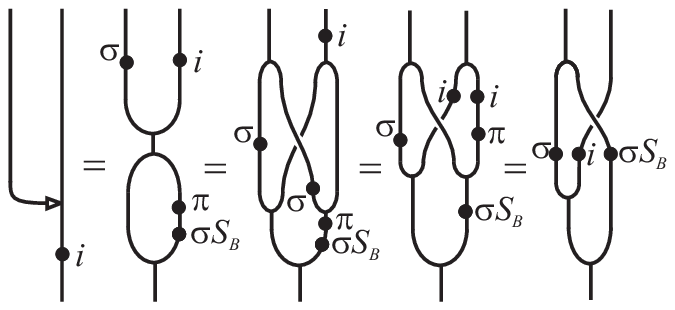}}
\end{center}
\caption{The proof of $i\tr=\mathrm{ad}$.}%
\label{fig:action}%
\end{figure}

c) We already know that $\tr$ is induced by the left adjoint
action. Obviously, if the right action $\tl$ is trivial then
$m_{R\Join B}=m_{R\# B}$, where $m_{R\# B}$ is defined in
(\ref{nn:smash}).
\end{proof}

\end{document}